\documentclass[12pt,oneside]{article}
\RequirePackage{epsfig}
\RequirePackage{amsmath}
\RequirePackage{amsfonts,amssymb}
\RequirePackage{amsthm}

\usepackage{color}
\usepackage{graphicx}
\usepackage{graphics}
\usepackage{subfigure}
\usepackage{footnote}
\usepackage{mathtools}
\newtheorem{lemma}{Lemma}
\newtheorem{proposition}{Proposition}
\newtheorem{defn}{Definition}
\newtheorem{theorem}{Theorem}
\newtheorem{rem}{Remark}
\newtheorem{coro}{Corollary}
\newtheorem{example}{Example}[section]

\newcommand{\ZZ}{\mathbb{Z}}
\newcommand{\RR}{\mathbb{R}}
\newcommand{\NN}{\mathbb{N}}
\newcommand{\CC}{\mathbb{C}}
\def\bdelta{\mbox{\boldmath $\delta$}}
\def\bgamma{\mbox{\boldmath $\gamma$}}
\def\bpi{\mbox{\boldmath $\pi$}}
\newcommand{\tp}{{\tilde p}}
\newcommand{\barr}{{\bar r}}
\newcommand\ba{\mathbf{a}}
\newcommand\bb{\mathbf{b}}
\newcommand\bbf{\mathbf{f}}

\newcommand\ie{{\it\thinspace i.e.}\ }

\begin{document}

\newpage

\begin{center}
\textbf{ON THE REPRODUCTION PROPERTIES OF NON-STATIONARY SUBDIVISION SCHEMES}

\bigskip

\textsc{Costanza Conti}\\
\footnotesize
\textrm{Universit\`{a} di Firenze, Dipartimento di Energetica ``Sergio Stecco''\\
Via Lombroso 6/17, 50134 Firenze, Italia
}\\
\textit{costanza.conti@unifi.it}

\bigskip
\normalsize
\textsc{Lucia Romani}\\
\footnotesize
\textrm{Universit\`{a} di Milano-Bicocca, Dipartimento di Matematica e Applicazioni\\
Via R. Cozzi 53, 20125 Milano, Italia}\\
\textit{lucia.romani@unimib.it}

\normalsize

\bigskip

\begin{abstract}
We present an accurate investigation of the algebraic conditions that the symbols of a convergent, univariate, binary, non-stationary subdivision scheme should fulfill in order to reproduce spaces of exponential polynomials.
A subdivision scheme is said to possess the property of reproducing exponential polynomials if, for any initial data uniformly sampled from some exponential polynomial function, the scheme yields the same function in the limit.
The importance of this property is due to the fact that several functions obtained as combinations of exponential polynomials (such as conic sections, spirals or special trigonometric and hyperbolic functions) are of great interest in graphical and engineering applications.
Since the space of exponential polynomials trivially includes standard polynomials, the results in this work extend the recently developed theory on polynomial reproduction to the non-stationary context. A significant application of the derived algebraic conditions on the subdivision symbols is the construction of new non-stationary subdivision schemes with specified reproduction properties.



\end{abstract}

\end{center}

\bigskip

\normalsize

{\renewcommand{\thefootnote}{}

\footnotetext{Date: February 17, 2010}
}

\section{Introduction}

Non-stationary subdivision schemes have proven to be efficient iterative algorithms to construct special classes of curves
ranging from polynomials and trigonometric curves to conic sections or spirals.
Aim of this paper is to establish the algebraic conditions that fully identify the reproduction properties of a given convergent, univariate, binary, non-stationary subdivision scheme.

\smallskip
Following the notation in \cite{DL02}, we denote by $\{\ba^{(k)},\ k\ge 0\}$ the finite sequence of real coefficients corresponding to the so called {\sl k-level mask} of the non-stationary subdivision scheme and we define by
\begin{equation}\label{symbol_def}
a^{(k)}(z)=\displaystyle{\sum_{j \in \ZZ} \, a^{(k)}_j z^j}, \quad k \geq 0, \quad z \in \CC \backslash\{0\}
\end{equation}
the Laurent polynomial whose coefficients are exactly the entries of $\ba^{(k)}$.
The polynomial in (\ref{symbol_def}) is commonly known as the {\sl k-level symbol}
of the non-stationary subdivision scheme.

\smallskip
Hereinafter we will denote through $\{S_{\ba^{(k)}},\ k\ge 0\}$ the linear subdivision operators based on the masks $\{\ba^{(k)},\ k\ge 0\}$, identifying the refinement process

\begin{equation}\label{def:suboper}
\bbf^{(k+1)}:=S_{\ba^{(k)}}\bbf^{(k)}, \quad (S_{\ba^{(k)}}\bbf^{(k)})_i:=\sum_{j\in \ZZ}a^{(k)}_{i-2j} \, f_j^{(k)}, \quad \ k\ge 0,
\end{equation}

\noindent
starting from any initial ``data" sequence $\bbf^{(0)} \equiv \bbf:=\{f_i \in \RR,\ i\in \ZZ\}$.
Note that, whenever the recursive relation in (\ref{def:suboper}) relies on the same
mask at each level of refinement - namely $\ba^{(k)}=\ba$ for all $k\geq 0$ - then the subdivision scheme is said stationary.

\smallskip
Attaching the data $f_i^{(k)}$ generated at the $k$-th step to the parameter values $t^{(k)}_i$ with
$$t^{(k)}_i<t_{i+1}^{(k)}, \quad \hbox{and} \quad t_{i+1}^{(k)}-t_i^{(k)}=2^{-k},\quad k\ge 0$$
we see that the subdivision process generates denser and denser sequences of data so that a notion of convergence can be established by
taking into account the piecewise
linear function $F^{(k)}$ that interpolates the data, namely
\[
  F^{(k)}(t_i^{(k)}) = f_i^{(k)}, \qquad
  F^{(k)}|_{[t_i^{(k)},t_{i+1}^{(k)}]} \in \Pi_1, \qquad
  i\in\ZZ,\quad k\geq0,
\]
where $\Pi_1$ is the space of linear polynomials. If
the sequence $\{F^{(k)},\ k\ge 0\}$ converges, then we denote its
limit by
\[
  g_\bbf := \lim_{k\to\infty} F^{(k)}
\]
and say that $g_\bbf$ is the \emph{limit function} of
the subdivision scheme based on the rule (\ref{def:suboper}) for the data $\bbf$. Due to the linearity of the refinement rule, we can equivalently check convergence of the subdivision scheme when applied to the initial data $\bdelta=\{\delta_{i,0}\ :i\in\ZZ\}=\{\dots,0,0,1,0,0,\dots\}$. If this is convergent with limit $\phi=g_{\bdelta}$,  then we have
\[
  g_\bbf = \sum_{j\in\ZZ} \phi(\cdot - j) f_j
\]
for any initial data sequence $\bbf$.

\smallskip Since most of the properties of a subdivision scheme (e.g. its convergence, its smoothness or its support size) do not depend on the choice of the parameter values $t_i^{(k)}$, these are usually set as
\begin{equation}\label{def:primal_par}
   t_i^{(k)}:=\frac{i}{2^k}, \quad i \in \ZZ, \quad k\ge 0.
\end{equation}
We refer to the choice in (\ref{def:primal_par}) as to the ``standard'' parametrization. As it will be better clarified later, with respect to the subdivision capability of reproducing specific classes of functions, the standard parametrization is not always the optimal one.
Indeed, the choice
\begin{equation}\label{def:dual_par}
   t_i^{(k)}:=\frac{i+p}{2^k},\quad i \in \ZZ, \quad p\in \RR,\quad k\ge 0
\end{equation}
with $p$ suitably set, turns out to be a better selection.
In particular, when $p=0$ we call this parametrization ``primal'', while in the case $p=-\frac{1}{2}$ we use the term ``dual'' parametrization, as in \cite{DHSS08}.
For a complete discussion concerning the choice of parametrization in the analysis of polynomial reproduction for stationary subdivision schemes, we refer the reader to the papers \cite{ContiHormann10, DHSS08}. In passing, we recall the following definition since it lays the foundations of our results in the non-stationary situation.

\begin{defn}\label{def:pg&pr}
A convergent stationary subdivision scheme $S_{{\ba}}$ is \emph{generating} polynomials up to degree $d_G$ if for any polynomial $\pi$ of degree $d\leq d_G$ there exists some initial data ${\bf q}^{(0)}$ such that $S^{\infty}_{\ba}{\bf q}^{(0)}=\pi$. Moreover, ${\bf q}^{(0)}$ is sampled from a polynomial of the same degree and with the same leading coefficient.
Additionally, a convergent subdivision scheme $S_{{\ba}}$ is \emph{reproducing} polynomials up to degree $d_R$ if for any polynomial $\pi$ of degree $d\leq d_R$ and for the initial data ${\bpi}^{(0)}=\{\pi(i),\ i\in\ZZ\}$ it results  $S^{\infty}_{\ba}{\bpi}^{(0)}=\pi$.
\end{defn}
We also recall that a convergent subdivision scheme that reproduces polynomials of degree $d_R$ has approximation order $d_R+1$ (see \cite{AdiLevin03} for the proof of this result).

\medskip \noindent In this context, the first purpose of this paper is to deal with the extension of the concepts of polynomial generation and polynomial reproduction to the non-stationary setting, leading to the concepts of exponential polynomial generation and exponential polynomial reproduction. In a few words, a subdivision scheme is said to possess the property of reproducing exponential polynomials if, for any initial data uniformly sampled from some exponential polynomial function, the scheme yields the same function in the limit. The second purpose of this paper is to show that, as in the stationary case discussed in \cite{ContiHormann10}, the choice of the correct parametrization is crucial in the non-stationary setting as well. \\
Our analysis brings to algebraic conditions, involving a parametrization of general type (\ref{def:dual_par}), that the symbols of a convergent, binary,  non-stationary subdivision scheme should fulfill in order to reproduce spaces of exponential polynomials. We remark that, as the symbol of the scheme changes from level to level and the
parametrization plays a crucial role in this kind of study, the proofs of the non-stationary case are often more difficult
than in the stationary case. This is why much of the results previously obtained can not be straightforwardly generalized but require a reformulation (see also \cite{DLL03}).\\
To illustrate the application of the simple but very general algebraic conditions derived, we construct novel subdivision symbols defining new non-stationary subdivision schemes with specified reproduction properties.

\medskip The outline of the paper  is the following. In Section \ref{sec:spaces}
we start by defining the space of exponential polynomials and we characterize some of its special instances that turn out to be really useful in applications.
Then, in Section \ref{sec:E-reproduction}, we provide new theoretical results concerning the algebraic conditions that the symbols of a convergent, binary, linear, non-stationary subdivision scheme should satisfy in order to reproduce functions from a given space of exponential polynomials. Finally, in Section \ref{sec:applications}, novel subdivision schemes are introduced to show how the illustrated conditions can be easily applied to work out non-stationary refinement rules with desired reproduction properties.

\section{Space of Exponential Polynomials}\label{sec:spaces}

\begin{defn}[Space of Exponential Polynomials] \label{def:Vspace}
Let $T\in \ZZ_+$ and $\bgamma=(\gamma_0,\gamma_1,\cdots,$ $\gamma_T)$
with $\gamma_T\neq 0$ a finite set of real or imaginary numbers and let $D^n$
the $n$-th order differentiation operator. The space of
exponential polynomials $V_{T,\bgamma}$ is the subspace
\begin{equation}\label{def:VMgamma}
V_{T,\bgamma}:=\{f:\RR\rightarrow \CC, f\in C^T(\RR):\quad
\sum_{j=0}^T \gamma_jD^j\,f =0\}.
\end{equation}
\end{defn}

\noindent A characterization of the space $V_{T,\bgamma}$ is provided by the following:

\begin{lemma}\cite{BR89} \label{Prop:characterization}
Let $\gamma(z)=\sum_{j=0}^T \gamma_jz^j$ and denote by
$\{\theta_\ell,\tau_\ell\}_{\ell=1,\cdots,N}$ the set of zeros
with multiplicity of $\gamma(z)$ satisfying
$$\gamma^{(r)}(\theta_\ell)=0,\quad
r=0,\cdots,\tau_\ell-1,\quad \ell=1,\cdots,N.$$ It results
$$
T=\sum_{\ell=1}^N \tau_\ell,\qquad
V_{T,\bgamma}:=Span\{x^{r}e^{\theta_\ell\,x},\
r=0,\cdots,\tau_\ell-1,\ \ \ell=1,\cdots,N\}.
$$
\end{lemma}

It is worth noting that given the space $V_{T,\bgamma}$ of
exponential polynomials,
we can define the space of exponential splines, say $SV_{T,\bgamma}$, as
the space of piecewise functions whose pieces are exponentials from the space $V_{T,\bgamma}$, opportunely joined with a certain degree of regularity.

\subsection{Characterization of special $V_{T,\bgamma}$ spaces}
In this section we illustrate some classical examples
of function spaces of the kind (\ref{def:VMgamma}).
In all the considered cases we assume $n\in \NN$ and we denote by $t$ a non-negative real or pure imaginary constant, namely
$$t\in \{0,s,{\rm i}s \, \vert s>0\}.$$

\noindent
\begin{example}\label{ex1}
{\rm
Let $N=3$ and $\theta_1=0,\ \theta_2=t,\ \theta_3=-t$ with multiplicities $\tau_1=n+1,\ \tau_2=\tau_3=1.$
It follows that $\{\theta_\ell,\tau_\ell\}_{\ell=1,2,3}$ identifies the set of zeros with multiplicities of the polynomial equation $\gamma(z)=z^{n+3}-t^2z^{n+1}=0$. Since the associated differential equation is $$D^{n+3}f(z)-t^2 D^{n+1}f(z)=0,$$
we have $$V_{n+3,\bgamma}=\{1,x,\cdots,x^{n},e^{tx},e^{-tx}\}.$$
The importance of this function space lies in the fact that it coincides with the following spaces:
\begin{itemize}
 \item {\emph{polynomials}} linear combination of $\{1,x,\cdots,x^{n},x^{n+1},x^{n+2}\}$, whenever $t=0$;
 \item {\emph{mixed polynomials-trigonometric functions}} linear combination of $\{1,x,\cdots,x^{n},$ $\cos(sx),\sin(sx)\}$, whenever $t={\rm i}s, \, s>0$;
\item {\emph{mixed polynomials-hyperbolic functions}} linear combination of $\{1,x,\cdots,x^{n},$ $\cosh(sx),\sinh(sx)\}$,
whenever $t=s>0$.
\end{itemize}
}
\end{example}

\begin{example}\label{ex2}
{\rm
We continue by introducing another interesting example of function space $V_{T,\bgamma}$.
Let $N=2n+1$ and $\theta_{\ell}, \ell=1,\cdots,2n+1$ the set of values $\theta_1=0,\quad \theta_{2j}=jt,\quad \theta_{2j+1}=-jt,\ j=1,\cdots,n$, with multiplicities $\tau_1=2,\quad \tau_{\ell}=1,\ \ell=2,\cdots,2n+1$.
They identify the set of zeros with multiplicities of the polynomial equation $\displaystyle{\gamma(z)=\sum_{i=0}^{2n+2}\gamma_iz^i=z^2\sum_{j=1}^n\left(z^2-j^2t^2\right)}$. Since the associated differential equation is $$\displaystyle{\sum_{i=0}^{2n+2}\gamma_i D^{i}f(z)=0},$$
we have
$$V_{2n+2,\bgamma}=\{1,x,e^{tx},e^{-tx}, \dots, e^{n tx}, e^{-n tx}\}.$$
The importance of this function space lies in the fact that it coincides with the space of
\begin{itemize}
 \item \emph{polynomials} linear combination of $\{1,x,\cdots,x^{2n},x^{2n+1}\}$, whenever $t=0$;
 \item \emph{mixed linear polynomials-trigonometric functions} linear combination of \\ $\{1,x,\cos(sx),\sin(sx),\cdots,\cos(n sx), \sin(n sx)\}$, whenever $t={\rm i}s, \, s>0$;
\item  \emph{mixed linear polynomials-hyperbolic functions} linear combination of \\ $\{1,x,\cosh(sx),\sinh(sx),\cdots,\cosh(n sx), \sinh(n sx)\}$, whenever $t=s>0$.
\end{itemize}
}
\end{example}

\begin{example}\label{ex3}
{\rm
We conclude with the illustration of a last important function space of type (\ref{def:VMgamma}). Let $N=3$, $\theta_1=0,\ \theta_2=t,\ \theta_3=-t$ with multiplicities $\tau_1=2,\ \tau_2=\tau_3=n.$
As $\{\theta_\ell,\tau_\ell\}_{\ell=1,2,3}$ identifies the set of zeros with multiplicities of the polynomial equation $\gamma(z)=z^2(z^2-t^2)^{n}=0$, associated with the differential equation $$\sum_{j=0}^{n}(-1)^j \, \binom{n}{j} \, t^{2j} D^{2n-2j+2}f(z)=0,$$ thus
$$V_{2n+2,\bgamma}=\{1,x,e^{tx},e^{-tx},xe^{tx},xe^{-tx},\cdots,x^{n-1}e^{tx},x^{n-1}e^{-tx}\}.$$
It coincides with the space of
\begin{itemize}
 \item \emph{polynomials} linear combination of $\{1,x,\cdots,x^{2n},x^{2n+1}\}$, whenever $t=0$;
 \item \emph{mixed linear polynomials-trigonometric functions} linear combination of \\ $\{1,x,\cos(sx),\sin(sx),\cdots,  x^{n-1}\cos(sx), x^{n-1}\sin(sx)\}$, whenever $t={\rm i}s, \, s>0$;
\item \emph{mixed linear polynomials-hyperbolic functions} linear combination of \\ $\{1,x,\cosh(sx),\sinh(sx),\cdots, x^{n-1}\cosh(sx), x^{n-1}\sinh(sx)\}$, whenever $t=s>0$.
\end{itemize}
}
\end{example}
\smallskip \noindent In practice, all the function spaces in the three examples can be used to reproduce conic sections.
Furthermore, the function space in Example \ref{ex2} has the additional advantage of reproducing trigonometric curves of higher order (such as the cardioid, the astroid, ...) while that in Example \ref{ex3} has the interesting capability of reproducing spirals. Obviously, for  all of them we can construct the corresponding $SV_{T,\bgamma}$ spaces.

\section{Exponential reproducing subdivision schemes}\label{sec:E-reproduction}

As it happens for polynomials, functions in the spaces $V_{T,\bgamma}$ can be generated via a subdivision recursion, but at the cost of considering non-stationary subdivision schemes. An important example of such non-stationary subdivision schemes
is defined by the following set of symbols
\begin{equation}\label{def:ES}
B_{n}^{(k)}(z)=2\prod_{i=1}^N\left(
\frac{e^{\frac{\theta_i}{2^{k+1}}}z+1}{e^{\frac{\theta_i}{2^{k+1}}}+1}\right)^{\tau_i},\quad
k\ge 0\,.
\end{equation}
As shown in \cite{MWW01, WW02}, the subdivision scheme with symbols in (\ref{def:ES}) generates limit functions belonging to the subclass of $C^{T-2}$ degree-$n$ L-splines (with $n=T-1$ and $T=\displaystyle{\sum_{i=1}^N \ \tau_i}$) whose pieces are exponentials
of the space $V_{T,\bgamma}$ where $\gamma(z)=\prod_{i=1}^N\left(
z-\theta_i\right)^{\tau_i}$. These functions, called
\emph{exponential B-splines}, are investigated in several papers (see, for instance \cite{MWW01, VonBluUnser07, WW02} and references therein).

\smallskip \noindent For later use, we note that the symbols in (\ref{def:ES}) for
$z^{(k)}_\ell:=e^{\frac{-\theta_\ell}{2^{k+1}}},\ \ell=1,\dots,N$, satisfy the conditions
$$
B_{n}^{(k)}(-z^{(k)}_\ell)=0,\quad \frac{d^r\, B_{n}^{(k)}(-z^{(k)}_\ell)}{dz^r}=0,\quad r=1,\dots,\tau_\ell-1.
$$
Furthermore, for $N \geq 1$ it turns out that
$$
B_{n}^{(k)}(z^{(k)}_\ell)=2\prod_{i=1}^N\left(
\frac{e^{\frac{\theta_i-\theta_\ell}{2^{k+1}}}+1}{e^{\frac{\theta_i}{2^{k+1}}}+1}\right)^{\tau_i}=2\left(
\frac{2}{e^{\frac{\theta_\ell}{2^{k+1}}}+1}\right)^{\tau_\ell}\prod_{i=1,\ i\neq \ell}^N\left(
\frac{e^{\frac{\theta_i-\theta_\ell}{2^{k+1}}}+1}{e^{\frac{\theta_i}{2^{k+1}}}+1}\right)^{\tau_i}\neq 2.
$$

\smallskip \noindent Notice that, when $\theta_{1}=0$ and $\tau_1=n+1$, then
$B_n^{(k)}(z)$ in (\ref{def:ES}) does not depend on $k$ and its
symbol is the stationary symbol of a degree-$n$ (polynomial) B-spline
$$
B_n(z)=\frac{(1+z)^{n+1}}{2^n}\,.
$$

\medskip \noindent We now continue with some definitions generalizing \cite[Definition 1.1, Definition 4.4.]{DHSS08} and successively we extend the conditions in \cite{ContiHormann10} to the binary non-stationary situation. Note that what we get is also an extension of \cite[Theorem 4.1]{ContiRomani09}, since the latter one is limited to the context of primal subdivision schemes only. We remind the reader that, within this section, the space $V_{T,\bgamma}$ is the one in Definition \ref{def:Vspace}.

\begin{defn} [E-Generation] \label{def:ERreproductionlimit}
Let $\{a^{(k)}(z),\ k\ge 0\}$ be a set of subdivision symbols. The subdivision scheme associated with the set of symbols $\{a^{(k)}(z),\ k\ge 0\}$ is said to be \emph{$V_{T,\bgamma}$-generating (for short E-generating)} if it is convergent and
for $f\in V_{T,\bgamma}$ and for the initial sequence ${\bbf}^{(0)}:=\{f(t^{(0)}_i),\ i\in \ZZ\}$,  it results
$$\lim_{k\rightarrow \infty}S_{\ba^{(k)}}\cdots S_{\ba^{(0)}}{\bbf}^{(0)}=\tilde f\,,\quad \tilde f\in V_{T,\bgamma}\,.
$$
\end{defn}

In the recent paper \cite{VonBluUnser07}, the authors proved a crucial result concerning E-generation under some mild assumptions on the subdivision symbols.

\begin{proposition}\cite[Theorem 1]{VonBluUnser07}\label{prop:ESgeneration}
Any convergent, non-stationary subdivision scheme associated with symbols $\{a^{(k)}(z),\ k\ge 0\}$ of finite length
that do not have pairs of opposite roots, is said to be generating a function in $V_{T,\bgamma}$ if and only if
$$
\frac{d^r\, a^{(k)}(-z^{(k)}_\ell)}{dz^r}=0,\quad r=0,\dots,\tau_\ell-1,\quad \hbox{for}\quad z^{(k)}_\ell:=e^{\frac{-\theta_\ell}{2^{k+1}}},\ \ell=1,\dots,N.
$$
\end{proposition}

\begin{rem}
As pointed out in \cite{VonBluUnser07}, Proposition \ref{prop:ESgeneration} states the important fact that a convergent, non-stationary subdivision scheme generates exponential polynomials if and only if its symbols $\{a^{(k)}(z),\ k\ge 0\}$ are divisible by the symbols of exponential B-splines $\{B_{n}^{(k)}(z),\ k\ge 0\}$ in (\ref{def:ES}). Moreover, exponential B-splines turn out to be the shortest functions satisfying these conditions.
\end{rem}

\noindent
Next we investigate the E-reproducing capability of a non-stationary subdivision scheme.
\smallskip

\begin{defn} [E-Reproduction] \label{def:ERreproductionlimit}
Let $\{a^{(k)}(z),\ k\ge 0\}$ be a set of subdivision symbols. The subdivision scheme associated with the symbols $\{a^{(k)}(z),\ k\ge 0\}$ is said to be \emph{$V_{T,\bgamma}$-reproducing (for short E-reproducing)} if it is convergent and
for $f\in V_{T,\bgamma}$ and for the initial sequence ${\bbf}^{(0)}:=\{f(t^{(0)}_i),\ i\in \ZZ\}$,  it results
$$\lim_{k\rightarrow \infty}S_{\ba^{(k)}}\cdots S_{\ba^{(0)}}{\bbf}^{(0)}=f\,.
$$
\end{defn}

\begin{defn} [Step-wise E-Reproduction] \label{def:ERreproduction_primal}
Let $\{a^{(k)}(z),\ k\ge 0\}$ be a set of subdivision symbols. The subdivision scheme associated with the symbols $\{a^{(k)}(z),\ k\ge 0\}$ is said to be \emph{step-wise $V_{T,\bgamma}$-reproducing} if
for $f\in V_{T,\bgamma}$ and for the sequences ${\bbf}^{(k+1)}:=\{f(t^{(k+1)}_i),\ i\in \ZZ\}$, $k\ge 0$,  it results
\begin{equation}\label{passo}
    \bbf^{(k+1)}=S_{\ba^{(k)}}\bbf^{(k)}\quad \hbox{or, equivalently,} \quad f (t^{(k+1)}_i)=\sum_{j\in \ZZ}a^{(k)}_{i-2j} \, f(t^{(k)}_j),\ \ i\in \ZZ\,.
\end{equation}
\end{defn}

\smallskip
\noindent
The next Proposition generalizes \cite[Corollary 4.5]{DHSS08} to the non-stationary situation. This is actually  already considered in \cite{ContiRomani09}, though in the context of ``primal'' subdivision schemes only.
However, since it is in fact independent of the parametrization, we just recall the statement and skip the proof.

\begin{proposition}\label{step-limit}\cite{ContiRomani09}
A convergent, non-stationary subdivision scheme based on the symbols $\{a^{(k)}(z),\ k\ge 0\}$ which is \emph{step-wise $V_{T,\bgamma}$-reproducing} is also \emph{$V_{T,\bgamma}$-reproducing} and vice versa.
\end{proposition}

\medskip \noindent We continue with a technical lemma that will be later used to prove the main result of this paper.

\begin{lemma}\label{lemma:polreprod_p}
Let  $\{a^{(k)}(z),\ k\ge 0\}$ be such that for $r=0,\dots,\tau_{\ell}-1$
\begin{equation}\label{eq:pol1}
  (z^{(k)}_\ell)^{p} \left( i+\frac p2 \right)^r=\sum_{j\in \ZZ}a^{(k)}_{2j}(i-j+p)^r\,(z^{(k)}_\ell)^{2j}\,,
\end{equation}
\begin{equation}\label{eq:pol2}
  (z^{(k)}_\ell)^{p}\left(i+\frac 12 +\frac p2 \right)^r=\sum_{j\in \ZZ}a^{(k)}_{2j+1}(i-j+p)^r\,(z^{(k)}_\ell)^{2j+1}\,.
\end{equation}
Then, for $r=0,\dots,\tau_{\ell}-1$
\begin{equation}\label{pari}
 \sum_{j\in \ZZ}a^{(k)}_{2j}(2j)^r\,(z^{(k)}_\ell)^{2j}=(p)^r (z^{(k)}_\ell)^{p}, \quad
 \sum_{j\in \ZZ}a^{(k)}_{2j}(2j+1)^r\,(z^{(k)}_\ell)^{2j}=(p+1)^r (z^{(k)}_\ell)^{p}
 \end{equation}
 and
 \begin{equation}\label{dispari}
 \sum_{j\in \ZZ}a^{(k)}_{2j+1}(2j)^r\,(z^{(k)}_\ell)^{2j+1}=(p-1)^r (z^{(k)}_\ell)^{p},\
 \sum_{j\in \ZZ}a^{(k)}_{2j+1}(2j+1)^r\,(z^{(k)}_\ell)^{2j+1}=(p)^r (z^{(k)}_\ell)^{p}.
 \end{equation}
\end{lemma}

\proof Expanding $(i+\frac p2)^r$ and $(i-j+p)^r$ in (\ref{eq:pol1}) we arrive at
\begin{equation}\label{pp}
 \left(z^{(k)}_\ell\right)^{p} \left(\frac p2 \right)^s=\sum_{j\in \ZZ}a^{(k)}_{2j}(p-j)^s\,(z^{(k)}_\ell)^{2j},\quad s=0,\dots,r,
\end{equation}
while expanding $(i+\frac12+\frac p2)^r$ and $(i-j+p)^r$ in (\ref{eq:pol2}) at
\begin{equation}\label{pd}
 \left(z^{(k)}_\ell\right)^{p} \left( \frac p2 +\frac1 2 \right)^s=\sum_{j\in \ZZ}a^{(k)}_{2j+1}(p-j)^s\,(z^{(k)}_\ell)^{2j+1},\quad s=0,\dots,r.
\end{equation}
Next, after substituting $\tp:=\frac p2$ in (\ref{pp}), we expand $(2\tp-j)^s$ and obtain
$$
\left(z^{(k)}_\ell\right)^{2\tp}\tp^s=\sum_{\barr=0}^s\left(
                                  \begin{array}{c}
                                    s \\
                                    \barr \\
                                  \end{array}
                                \right)
 (2\tp)^{s-\barr} \left( -\frac12 \right)^\barr\sum_{j\in \ZZ}a^{(k)}_{2j}(2j)^\barr\,(z^{(k)}_\ell)^{2j}, \,
\begin{array}{l}
 s=0,\dots,r, \\ r=0,\dots,\tau_{\ell}-1.
\end{array}
$$
Now, since the above equality is possible if and only if
$$
\left( -\frac12 \right)^\barr\sum_{j\in \ZZ}a^{(k)}_{2j}(2j)^\barr\,(z^{(k)}_\ell)^{2j}=(-\tp)^\barr \left(z^{(k)}_\ell\right)^{2\tp},
$$
we arrive at the first condition in (\ref{pari}) \ie,
$$\left(z^{(k)}_\ell\right)^{p}(p)^r=\sum_{j\in \ZZ}a^{(k)}_{2j}(2j)^r\,(z^{(k)}_\ell)^{2j},\quad r=0,\dots,\tau_{\ell}-1.$$
Next, expanding $(2j+1)^r$ and using the first condition in (\ref{pari}) we obtain the second condition in (\ref{pari}) for $r=0,\dots,\tau_{\ell}-1$:
$$
\sum_{j\in \ZZ}a^{(k)}_{2j}(2j+1)^r\,(z^{(k)}_\ell)^{2j}=\sum_{\barr=0}^r\left(
                                  \begin{array}{c}
                                    r \\
                                    \barr \\
                                  \end{array}
                                \right)\sum_{j\in \ZZ}a^{(k)}_{2j}(2j)^\barr\,(z^{(k)}_\ell)^{2j}=\left(z^{(k)}_\ell\right)^{p}(p+1)^r.
$$
\smallskip
To prove the second part of the Lemma, with the substitution $\tp:=\frac p2 $ we write (\ref{pd}) as
\begin{equation}\label{pd2}
\left(z^{(k)}_\ell\right)^{2\tp} \left( \tp+\frac 12 \right)^s=\sum_{j\in \ZZ}a^{(k)}_{2j+1}(2\tp-j)^s\,(z^{(k)}_\ell)^{2j+1},\quad s=0,\dots,r.
\end{equation}
\smallskip
Then, if we expand $(2\tp-j)^s$ we arrive at
$$
\left(z^{(k)}_\ell\right)^{2\tp} \left( \tp+\frac{1}{2} \right)^s=\sum_{\barr=0}^s\left(
                                  \begin{array}{c}
                                    s \\
                                    \barr \\
                                  \end{array}
                                \right)
 (2\tp)^{s-\barr} \left(-\frac{1}{2} \right)^\barr \sum_{j\in \ZZ}a^{(k)}_{2j+1}(2j)^\barr\,(z^{(k)}_\ell)^{2j+1}
$$
and therefore, since the above equality is possible if and only if
$$
\left( -\frac12 \right)^\barr\sum_{j\in \ZZ}a^{(k)}_{2j+1}(2j)^\barr\,(z^{(k)}_\ell)^{2j+1}=\left(-\tp+\frac12 \right)^\barr \left(z^{(k)}_\ell\right)^{2\tp},
$$  we get
$$
 \left(z^{(k)}_\ell\right)^{2\tp} (2\tp-1)^\barr=\sum_{j\in \ZZ}a^{(k)}_{2j+1}(2j)^\barr\,(z^{(k)}_\ell)^{2j+1}, \quad \barr=0,\dots,s,\ s=0,\dots,\tau_{\ell}-1
$$
which coincides with the first condition in (\ref{dispari}):
$$
 \left(z^{(k)}_\ell\right)^{p} (p-1)^r=\sum_{j\in \ZZ}a^{(k)}_{2j+1}(2j)^r\,(z^{(k)}_\ell)^{2j+1},\quad r=0,\dots, \tau_{\ell}-1.
$$
Last, expanding $(2j+1)^r$ and using the just proven first condition in (\ref{dispari}), we obtain the second condition in (\ref{dispari}) for $r=0,\dots, \tau_{\ell}-1$
$$
\sum_{j\in \ZZ}a^{(k)}_{2j+1}(2j+1)^r \,(z^{(k)}_\ell)^{2j+1}=\sum_{\barr=0}^r\left(
                                  \begin{array}{c}
                                    r \\
                                    \barr \\
                                  \end{array}
                                \right)\sum_{j\in \ZZ}a^{(k)}_{2j+1}(2j)^\barr\,(z^{(k)}_\ell)^{2j+1}=(p)^r (z_{\ell}^{(k)})^p\,.
$$
\qed

\noindent
We proceed with another useful intermediate result.

\begin{proposition}\label{prop:coeff_der}
For the polynomial coefficients $\{p_{r,s},\ d_{r,s}\}_{s=0, \cdots,r}$ satisfying

\begin{equation}\label{coeff-der1}
\sum_{s=0}^rp_{r,s}(2j)^s:=\prod_{q=0}^{r-1} (2j-q)\quad \hbox{and}\quad
\sum_{s=0}^rd_{r,s}(2j)^s:=\prod_{q=0}^{r-1} (2j+1-q),
\end{equation}
it holds
\begin{equation}\label{coeff-der2}
\sum_{s=0}^rp_{r,s}(p)^s=\prod_{i=0}^{r-1}(p-i)\quad \hbox{and}\quad \sum_{s=0}^rd_{r,s}(p-1)^s=\prod_{i=0}^{r-1}(p-i)\,.
\end{equation}
\end{proposition}

\proof The proof is trivial since the first equality simply follows by replacing $(2j)$ with $p$ in the first  of (\ref{coeff-der1}), while the second one by replacing $(2j)$ with $p-1$. \qed

\bigskip \noindent
We are now in a position to state the central result of this paper, which is an extension of the main theorem in \cite{ContiHormann10} to the binary, non-stationary context.

\begin{theorem}\label{theo:ERreproduction}
Let $z^{(k)}_\ell:=e^{\frac{-\theta_\ell}{2^{k+1}}},\ \ell=1,\cdots,N$. A convergent subdivision scheme associated with  $\{a^{(k)}(z),\ k\ge 0\}$ reproduces $V_{T,\bgamma}$ if and only if for each $k\ge 0$ it holds
\begin{equation}\label{ERcond}
\begin{array}{ll}
      a^{(k)}(z^{(k)}_\ell)=2 \, \left(z^{(k)}_{\ell}\right)^p, \quad  \frac{d^r\, a^{(k)}(z^{(k)}_\ell)}{dz^r}=2\left(z^{(k)}_\ell\right)^{p-r}\,\displaystyle{\prod_{i=0}^{r-1}(p-i)},\quad
      \begin{array}{l}
      \ell=1,\dots,N,\\
      r=1,\dots,\tau_\ell-1,
      \end{array}
      \\
      \\
      a^{(k)}(-z^{(k)}_\ell)=0,\quad  \frac{d^r\, a^{(k)}(-z^{(k)}_\ell)}{dz^r}=0,\quad
      \begin{array}{l}
      \ell=1,\dots,N,\\
      r=1,\dots,\tau_\ell-1.
      \end{array}
    \end{array}
\end{equation}
\end{theorem}
\proof
Let us consider the function $f(x)=x^r\, e^{\theta_{\ell}\, x},\ r\in \{0,\dots, \tau_{\ell}-1\}$ with $\theta_{\ell}$ a root of multiplicity $\tau_{\ell}$ of $\gamma(z)$.
It obviously holds that $f\in V_{T,\bgamma}$. Let $\{a^{(k)}(z),\ k\ge 0\}$ be an E-reproducing non-stationary subdivision scheme. Then, from (\ref{passo}) it holds
(splitting the even and odd elements)
$$
\left( \frac {2i+p}{2^{k+1}} \right)^r\,e^{\theta_{\ell}\, \frac {2i+p}{2^{k+1}}}=\sum_{j\in \ZZ}a^{(k)}_{2j} \left( \frac {i-j+p}{2^{k}} \right)^r\, e^{\theta_{\ell} \, \frac {i-j+p}{2^{k}}}\,,
$$
$$
\left( \frac {2i+1+p}{2^{k+1}} \right)^r\,e^{\theta_{\ell} \, \frac {2i+1+p}{2^{k+1}}}=\sum_{j\in \ZZ}a^{(k)}_{2j+1} \left(\frac {i-j+p}{2^{k}} \right)^r\, e^{\theta_{\ell} \, \frac {i-j+p}{2^{k}}}\,,
$$
or, equivalently
\begin{equation}\label{eq:polripetizione1}
 \left( i+\frac p2 \right)^r=\underbrace{e^{\theta_{\ell}\, \frac {p}{2^{k+1}}}}_{\left(z^{(k)}_\ell\right)^{-p}}\sum_{j\in \ZZ}a^{(k)}_{2j}(i-j+p)^r\,(z^{(k)}_\ell)^{2j}\,,
\end{equation}
\begin{equation}\label{eq:polripetizione2}
 \left( i+\frac 12 +\frac p2 \right)^r=\underbrace{e^{\theta_{\ell}\, \frac {p}{2^{k+1}}}}_{\left(z^{(k)}_\ell\right)^{-p}}\sum_{j\in \ZZ}a^{(k)}_{2j+1}(i-j+p)^r\,(z^{(k)}_\ell)^{2j+1}\,.
\end{equation}
Taking (\ref{eq:polripetizione1},\ref{eq:polripetizione2}) for $r=0$ and summing them up, we obtain
$$
\sum_{j\in \ZZ}a^{(k)}_{j}(z^{(k)}_\ell)^{j}=\sum_{j\in \ZZ}a^{(k)}_{2j}(z^{(k)}_\ell)^{2j}+\sum_{j\in \ZZ}a^{(k)}_{2j+1}(z^{(k)}_\ell)^{2j+1}=2\,\left(z^{(k)}_\ell\right)^{p}\,,
$$
as well as
$$
\sum_{j\in \ZZ}a^{(k)}_{j}(-z^{(k)}_\ell)^{j}=\sum_{j\in \ZZ}a^{(k)}_{2j}(-z^{(k)}_\ell)^{2j}+\sum_{j\in \ZZ}a^{(k)}_{2j+1}(-z^{(k)}_\ell)^{2j+1}=0\,,
$$
which are the first relations in (\ref{ERcond}).
To prove the other relations in (\ref{ERcond}) we proceed in the following way.
For $r \in \{1, ..., \tau_{\ell}-1\}$, taking the $r-th$ derivative of $$a^{(k)}(z)=\sum_{j\in \ZZ}a^{(k)}_{2j}z^{2j}+\sum_{j\in \ZZ}a^{(k)}_{2j+1}z^{2j+1}$$
after multiplication  by $z^r$ we have
$$
z^r\frac{d^r\, a^{(k)}(z)}{dz^r}=\sum_{j\in \ZZ}a^{(k)}_{2j}\left(\prod_{q=0}^{r-1} (2j-q) \right)z^{2j}+\sum_{j\in \ZZ}a^{(k)}_{2j+1}\left(\prod_{q=0}^{r-1} (2j+1-q) \right)z^{2j+1}.
$$
Therefore, using Lemma \ref{lemma:polreprod_p} and Proposition \ref{prop:coeff_der} for several times, we obtain

\begin{displaymath}
\begin{array}{ll}
\left(z^{(k)}_\ell\right)^{r}\frac{d^{r}\, a^{(k)}(z^{(k)}_\ell)}{dz^{r}}&=\displaystyle{\sum_{j\in \ZZ}a^{(k)}_{2j}[2j-(r-1)]\left(\prod_{q=0}^{r-2} (2j-q) \right)(z^{(k)}_\ell)^{2j}}
\\ \\
& +\displaystyle{\sum_{j\in \ZZ}a^{(k)}_{2j+1}[2j+1-(r-1)]} \left(\prod_{q=0}^{r-2} (2j+1-q) \right)(z^{(k)}_\ell)^{2j+1}\\
\\
&=\displaystyle{\sum_{s=0}^{r-1} p_{r-1,s}\sum_{j\in \ZZ}a^{(k)}_{2j}(2j)^{s+1}(z_{\ell}^{(k)})^{2j}}\\
\\
&+\displaystyle{\sum_{s=0}^{r-1} d_{r-1,s}\sum_{j\in \ZZ}a^{(k)}_{2j+1}(2j)^{s+1}(z_{\ell}^{(k)})^{2j+1}}\\
\\
&-(r-1)\displaystyle{\sum_{s=0}^{r-1} p_{r-1,s}\sum_{j\in \ZZ}a^{(k)}_{2j}(2j)^{s}(z_{\ell}^{(k)})^{2j}}
\\\\
\end{array}
\end{displaymath}

\begin{equation}\label{der_r}
\begin{array}{ll}
&+[1-(r-1)]\displaystyle{\sum_{s=0}^{r-1} d_{r-1,s}\sum_{j\in \ZZ}a^{(k)}_{2j+1}(2j)^{s}(z^{(k)}_\ell)^{2j+1}}\\
\\
&=\left(z^{(k)}_\ell\right)^{p}\left(\displaystyle{p\sum_{s=0}^{r-1} p_{r-1,s}(p)^{s}}+(p-1)\displaystyle{\sum_{s=0}^{r-1} d_{r-1,s}(p-1)^{s}}-(r-1)\,\displaystyle{\sum_{s=0}^{r-1} p_{r-1,s}(p)^{s}} \right.
\\ \\
& + \left. [1-(r-1)]\displaystyle{\sum_{s=0}^{r-1} d_{r-1,s}(p-1)^{s}}\right)
\\ \\
&=\left(z^{(k)}_\ell\right)^{p}[p-(r-1)]\,\left(\displaystyle{\sum_{s=0}^{r-1} p_{r-1,s}(p)^{s}}+\displaystyle{\sum_{s=0}^{r-1} d_{r-1,s}(p-1)^{s}}\right)\\
\\
&=2 \left(z^{(k)}_\ell\right)^{p}\,\displaystyle{\prod_{i=0}^{r-1}(p-i)}\,.
\\
\end{array}
\end{equation}
\normalsize
Hence,
$$
\frac{d^{r}\, a^{(k)}(z^{(k)}_\ell)}{dz^{r}}=2\left(z^{(k)}_\ell\right)^{p-r}\,\displaystyle{\prod_{i=0}^{r-1}(p-i)}\,.
$$
Differently, by replacing $z^{(k)}_\ell$ with $-z^{(k)}_\ell$ in (\ref{der_r}), we obtain

$$
\small
\begin{array}{ll}
\left(-z^{(k)}_\ell\right)^{r}\frac{d^{r}\, a^{(k)}(-z^{(k)}_\ell)}{dz^{r}}&=\left(z^{(k)}_\ell\right)^{p}\left(\displaystyle{p\sum_{s=0}^{r-1} p_{r-1,s}(p)^{s}}-(p-1)\displaystyle{\sum_{s=0}^{r-1} d_{r-1,s}(p-1)^{s}}+ \right.
\\ \\
& \left. -(r-1)\,\displaystyle{\sum_{s=0}^{r-1} p_{r-1,s}(p)^{s}}-[1-(r-1)]\displaystyle{\sum_{s=0}^{r-1} d_{r-1,s}(p-1)^{s}}\right)\\
\\
&=\left(z^{(k)}_\ell\right)^{p}[p-(r-1)]\,\left(\displaystyle{\sum_{s=0}^{r-1} p_{r-1,s}(p)^{s}}-\displaystyle{\sum_{s=0}^{r-1} d_{r-1,s}(p-1)^{s}}\right)\\
\\
&=0\,,
\\
\end{array}
$$
\normalsize
which concludes the proof of the necessary part of the claim.

\smallskip \noindent To prove the converse direction similar arguments can be used.\qed

\medskip
\begin{rem}
An important consequence of the previous Theorem is that it allows us to identify, whenever it exists, the correct parametrization for a
given non-stationary subdivision scheme. In fact $p$ has to satisfy
\begin{equation}\label{cond:correctp}
    a^{(k)}(z^{(k)}_\ell)=2(z_\ell^{(k)})^p,\quad \ell=1,\dots,N.
\end{equation}
On the other hand, in case for some $\tilde \ell$ it results $z_{\tilde \ell}^{(k)}=1$ and $\tau_{\ell}>1$ for all $\ell=1,\cdots, N$, by fixing $r=1$ we get simpler conditions to be satisfied by $p$:
\begin{equation}\label{cond:correctp}
  p=\frac{1}{2}\frac{d a^{(k)}(z^{(k)}_{\tilde \ell})}{dz}\quad \hbox{and}\quad
   p=\frac{1}{2\left(z^{(k)}_\ell\right)^{p-1}}\frac{d a^{(k)}(z^{(k)}_\ell)}{dz},\quad \ell=1,\dots,N, \quad \ell\neq \tilde \ell.
\end{equation}
\end{rem}

\begin{defn} [Symmetric Schemes] \label{def_symm}
A $k$-level subdivision mask ${\bf a}^{(k)}$ is called {\rm odd symmetric} if $a^{(k)}_{-i}=a^{(k)}_{i}$, $i \in \ZZ$ and {\rm even symmetric} if $a^{(k)}_{-i}=a^{(k)}_{i-1}$, $i \in \ZZ$. In terms of $k$-level symbol these conditions read as $a^{(k)}(z)=a^{(k)}(z^{-1})$ and $a^{(k)}(z)z=a^{(k)}(z^{-1})$, respectively.
\end{defn}

\medskip \noindent
From Theorem \ref{theo:ERreproduction} the following result follows straightforwardly (see also \cite{DHSS08}).

\begin{coro}
If a convergent, non-stationary subdivision scheme with symbols $\{a^{(k)}(z),$ $ k\geq 0\}$ satisfies conditions (\ref{ERcond}) for the root $z_1^{(k)}=1$ (namely reproduces constants) and it is odd (respectively, even) symmetric, then it reproduces functions from the space of exponential polynomials $V_{T,\bgamma}$ with respect to the primal (respectively, dual) parametrization.
\end{coro}

\medskip
We further observe that Theorem \ref{theo:ERreproduction} applies to interpolatory symbols with the parametrization choice $p=0$ and thus, in this special case, it gives back the results in \cite[Theorem 4.1]{ContiRomani09} and in \cite[Theorem 2.3]{DLL03}.

\begin{coro}
An interpolatory scheme reproduces functions from the space of exponential polynomials $V_{T,\bgamma}$ with respect to the primal parametrization $p=0$.
\end{coro}

\proof
Since interpolatory symbols satisfy the relation $a^{(k)}(z)+a^{(k)}(-z)=2$, by taking the first derivative on both sides of this relation it turns out that
\begin{equation}\label{cond_der_interp}
\frac{d}{dz}a^{(k)}(z)-\frac{d}{dz}a^{(k)}(-z)=0.
\end{equation}
Hence by evaluating (\ref{cond_der_interp}) at an arbitrary $z_{\ell}^{(k)}:=e^{-\frac{\theta_{\ell}}{2^{k+1}}}$, it follows that $\frac{d}{dz}a^{(k)}(z_{\ell}^{(k)})=\frac{d}{dz}a^{(k)}(-z_{\ell}^{(k)})$. Taking into account conditions (\ref{ERcond}) this holds true only in the case $p=0$.
\qed

\bigskip
We close this section by showing that a shift of $z^n$ in all the symbols $\{a^{(k)}(z), \ k \geq 0\}$ of an E-reproducing subdivision scheme based on the parametrization $p$, yields the replacement $p \rightarrow p+n$ in conditions (\ref{ERcond}),  that is an E-reproducing subdivision scheme based on the parametrization $p+n$. See \cite{ContiHormann10} for the analogous result in the stationary case.

\begin{proposition}\label{cor:shift2}
If $S_{\ba^{(k)}},\ k\ge 0$ is a subdivision scheme with symbols $\{a^{(k)}(z), \ k \geq 0\}$ that reproduces a function space $V_{T, \bgamma}$, then so does the scheme $S_{\bb^{(k)}},\ k\ge 0$ with symbols $b^{(k)}(z)=a^{(k)}(z)z^n,\ k\ge 0$ for any $n\in\ZZ$.
\end{proposition}

\proof
According to the Leibnitz rule, the $r$-th derivative of the symbol $b^{(k)}(z)$ is
\begin{equation}\label{eq:b-derivative}
 \frac{d^{r}\,  b^{(k)}(z)}{dz^{r}} =
    \sum_{i=0}^r \binom{r}{i} \frac{d^{i}a^{(k)}(z)}{dz^{i}} \prod_{j=0}^{r-i-1} (n-j) z^{n-r+i}.
\end{equation}
For $z=-z_\ell^{(k)}$ and $\ell= 1,\dots,N$, each addend of the sum in (\ref{eq:b-derivative}) vanishes because $\frac{d^{i}a^{(k)}(-z_\ell^{(k)})}{dz^{i}}=0$ for $i=0,\dots,\tau_\ell-1,\ \ell=1,\cdots,N$ by assumption, and so
$$\frac{d^{r}b^{(k)}(-z_\ell^{(k)})}{dz^{r}}=0, \quad r=0, \cdots, \tau_\ell-1, \quad \ell=1,\cdots,N.$$
To complete the proof, we show by induction that
\begin{equation}\label{eq:b-reproduction}
  \frac{d^{r}b^{(k)}(z_\ell^{(k)})}{dz^{r}}= 2\left( z_\ell^{(k)} \right)^{p+n-r}\prod_{i=0}^{r-1}(p+n-i),
  \qquad r=0,\dots,\tau_\ell-1,
\end{equation}
where $p$ is the parametric shift of $S_{\ba^{(k)}}$. For $r=0$ this holds because $b^{(k)}(z_\ell^{(k)})=a^{(k)}(z_\ell^{(k)})\left(z_\ell^{(k)}\right)^n=2\left(z_\ell^{(k)}\right)^{p+n}$. Now let $0 < r \leq \tau_\ell-1$.
Since when $z=z_\ell^{(k)}$ the first conditions in (\ref{ERcond}) can be rewritten in a recursive way as
$$\frac{d^{i+1}\, a^{(k)}(z^{(k)}_\ell)}{dz^{i+1}}=\frac{(p-i)}{z_\ell^{(k)}} \frac{d^{i}\, a^{(k)}(z^{(k)}_\ell)}{dz^{i}},\quad \ell=1,\dots,N,$$ exploiting the fact that $\binom{r+1}{i}=\binom{r}{i}+\binom{r}{i-1}$, we find that

\begin{align*}
 \frac{d^{r+1}b^{(k)}(z_\ell^{(k)})}{dz^{r+1}}
  &= \sum_{i=0}^{r+1} \binom{r+1}{i} \frac{d^{i}a^{(k)}(z_\ell^{(k)})}{dz^{i}} \prod_{j=0}^{r-i} (n-j)\left(z_\ell^{(k)}\right)^{n-(r+1)+i}\\
  \\
  &= \sum_{i=0}^{r+1} \binom{r}{i} \frac{d^{i}a^{(k)}(z_\ell^{(k)})}{dz^{i}} \prod_{j=0}^{r-i} (n-j)\left(z_\ell^{(k)}\right)^{n-(r+1)+i}\\
  \\
  & + \sum_{i=0}^{r+1} \binom{r}{i-1} \frac{d^{i}a^{(k)}(z_\ell^{(k)})}{dz^{i}} \prod_{j=0}^{r-i} (n-j)\left(z_\ell^{(k)}\right)^{n-(r+1)+i}\\
   \\
  & = \sum_{i=0}^{r} \binom{r}{i} \frac{d^{i}a^{(k)}(z_\ell^{(k)})}{dz^{i}} \prod_{j=0}^{r-i} (n-j)\left(z_\ell^{(k)}\right)^{n-(r+1)+i}\\
   \\
   & + \sum_{i=0}^{r} \binom{r}{i} \frac{d^{i+1}a^{(k)}(z_\ell^{(k)})}{dz^{i+1}}\prod_{j=0}^{r-i-1} (n-j)\left(z_\ell^{(k)}\right)^{n-(r+1)+i+1}\\
     \\
  \end{align*}
  \begin{align*}
     &= \sum_{i=0}^{r} \binom{r}{i} \frac{d^{i}a^{(k)}(z_\ell^{(k)})}{dz^{i}} \prod_{j=0}^{r-i} (n-j)\left(z_\ell^{(k)}\right)^{n-(r+1)-i}\\
  \\
    & + \sum_{i=0}^{r} \binom{r}{i} \frac{d^{i}a^{(k)}(z_\ell^{(k)})}{dz^{i}} \frac{(p-i)}{z_\ell^{(k)}} \prod_{j=0}^{r-i-1} (n-j)\left(z_\ell^{(k)}\right)^{n-(r+1)-i+1}\\
   \\
  &= \sum_{i=0}^{r} \binom{r}{i} \frac{d^{i}a^{(k)}(z_\ell^{(k)})}{dz^{i}} \left(z_\ell^{(k)}\right)^{n-(r+1)-i}\prod_{j=0}^{r-i-1} (n-j)
     \bigl( (n-r+i) + (p-i) \bigr)\\
  &=\frac{(p+n-r)}{z_\ell^{(k)}}\, \frac{d^{r}b^{(k)}(z_\ell^{(k)})}{dz^{r}}.
\end{align*}

\noindent
Hence, if \eqref{eq:b-reproduction} is true for $r$, then it is also true for $r+1$, so concluding the proof.
\qed

%
\section{Application examples}\label{sec:applications}
In this section we consider some examples of new E-generating and E-reproducing non-stationary subdivision schemes. The latter ones are obtained by the application of Theorem \ref{theo:ERreproduction} which is also exploited to compute the correct parametrization to be associated with these schemes.

\begin{example} {Primal scheme generating $V_{T,\bgamma}=\{1, \ x, \ e^{tx},\ e^{-tx}\}$.}

{\rm
\bigskip
\noindent
We start by considering the set of symbols $\{a_1^{(k)}(z),\ k \geq 0\}$ with
\begin{equation}\label{newscheme1}
\begin{array}{l}
a_1^{(k)}(z)=z^{-3} \, \frac{(z+1)^2}{2} \, \frac{z^2+2v^{(k)}z+1}{2(v^{(k)}+1)} \, \frac{\left(2+\sqrt{2(v^{(k)}+1)}\right)z^2+2\left(2(v^{(k)}+2)+3\sqrt{2(v^{(k)}+1)}\right)z+\left(2+\sqrt{2(v^{(k)}+1)}\right)}{4\left(v^{(k)}+3+2\sqrt{2(v^{(k)}+1)}\right)}
\end{array}
\end{equation}
where, for all $k \geq 0$, $v^{(k)}$ is updated through the recurrence formula
\begin{equation}\label{rec_par}
v^{(k)}=\sqrt{\frac{1+v^{(k-1)}}{2}}\quad \hbox{with} \quad v^{(-1)}>-1.
\end{equation}
To investigate the generation/reproduction properties of this scheme we use Proposition \ref{prop:ESgeneration} and Theorem \ref{theo:ERreproduction} referring to the function space $V_{T,\bgamma}=\{1,\ x, \ e^{tx},\ e^{-tx}\}$ corresponding to the values $\theta_1=0,\ \theta_2=t,\theta_3=-t$ and  $\tau_1=2,\ \tau_2=1,\ \tau_3=1$ .
First of all we verify that $a_1^{(k)}(z)$ does not have pairs of opposite roots.
Then, with some algebra we easily see that for $$z^{(k)}_1=1,\quad z^{(k)}_2=e^{t/2^{k+1}},\quad z^{(k)}_3=e^{-t/2^{k+1}},\quad t\in \{0,s,\mathrm{i}s\,|\, s>0\}$$
the conditions
$$
a_1^{(k)}(-z^{(k)}_\ell) =0,\ \ell=1,\, 2,\, 3,
$$
and the condition
$$
\frac{d}{dz}a_1^{(k)}(-z_1^{(k)})=0,
$$
are satisfied so that the scheme, indeed, generates the space $V_{T,\bgamma}$.
However, such a scheme does not reproduce the whole space above but only the subspace $\{1, \ x\}$ since conditions in Theorem \ref{theo:ERreproduction} concerning the roots $z_\ell^{(k)}$, $\ell=2,3$ are not fulfilled.\\
Finally, we observe that such a scheme is $C^3$ continuous since asymptotically equivalent \cite{DL95} to the stationary $C^3$ scheme with symbol $$\frac{(z+1)^4}{2^3} \left(\frac{1}{8}z^2 + \frac{3}{4}z + \frac{1}{8}\right)=\frac{(z+1)^3}{2^3} (z+1)\left(\frac{1}{8}z^2 + \frac{3}{4}z + \frac{1}{8}\right).$$
}
\end{example}

\medskip
The next example is to provide a primal scheme reproducing the entire space $V_{T,\bgamma}=\{1, \ x, \ e^{tx},\ e^{-tx}\}$. To this purpose, to keep the E-generation of the space $V_{T,\bgamma}$,  we multiply the symbol in (\ref{newscheme1}) by the factor $\alpha z^{-1}+(1-2\alpha)+\alpha z$ and determine the value of $\alpha$ such that the required conditions on $z_2^{(k)}$, $z_3^{(k)}$ are also verified.
\medskip

\begin{example} {Primal scheme reproducing $V_{T,\bgamma}=\{1, \ x, \ e^{tx},\ e^{-tx}\}$.}

{\rm
\bigskip
\noindent
Let

$$
a_2^{(k)}(z)=a_1^{(k)}(z) \cdot (\alpha z^{-1}+(1-2\alpha)+\alpha z), \quad {\rm with} \quad
\alpha=\frac{2-v^{(k)} \sqrt{2(v^{(k)}+1)}}{2v^{(k)} (v^{(k)} - 1)\sqrt{2(v^{(k)}+1)}}.
$$
\begin{equation}\label{newscheme2}
\end{equation}

\medskip
\noindent
Clearly, due to the factor $a_1^{(k)}(z)$, the scheme with $k$-level symbol $a_2^{(k)}(z)$ generates the function space $V_{T,\bgamma}=\{1, \ x, \ e^{tx},\ e^{-tx}\}$. Our intention is now to show that, differently from $a_1^{(k)}(z)$, the latter one exactly reproduces that space.
First of all, by exploiting conditions (\ref{cond:correctp}) we identify that the correct parametrization is given by $p=0$
and hence the scheme is primal.
Further, since conditions
$$
a_2^{(k)}(-z^{(k)}_\ell) =0, \quad a_2^{(k)}(z^{(k)}_\ell)=2, \quad \ell= 1, \ 2,\, 3,
$$
and
$$
\frac{d}{dz}a_2^{(k)}(-z_1^{(k)})=\frac{d}{dz}a_2^{(k)}(z_1^{(k)})=0,
$$
are satisfied, the scheme indeed reproduces the space $V_{T,\bgamma}$ with respect to the primal parametrization $p=0$.

We conclude the example by noticing  that the non-stationary subdivision scheme (\ref{newscheme2}) is of class $C^{2}$ since asymptotically equivalent \cite{DL95} to the stationary $C^2$ scheme with symbol

\begin{displaymath}
\begin{array}{ll}
& \frac{(z + 1)^4}{2^3} \left(-\frac{5}{64}z^4-\frac{3}{16}z^3+\frac{49}{32}z^2-\frac{3}{16}z-\frac{5}{64} \right)=\\ \\
& \frac{(z + 1)^2}{2^2} (z+1)\left(-\frac 5{128} z^5 -\frac {17}{128} z^4+\frac{43}{64}z^3+\frac{43}{64}z^2 -\frac {17}{128} z -\frac 5{128} \right).
\end{array}
\end{displaymath}
}
\end{example}

\medskip
Aim of the forthcoming example is to define a dual scheme that can reproduce the same function space $V_{T,\bgamma}=\{1, \ x, \ e^{tx},\ e^{-tx}\}$.
\medskip

\begin{example} {Dual scheme reproducing $V_{T,\bgamma}=\{1, \ x, \ e^{tx},\ e^{-tx}\}$.}

{\rm
\bigskip
\noindent
We consider the $k$-level symbol

\begin{equation}\label{newscheme3}
\begin{array}{ll}
& a_3^{(k)}(z)=z^{-4} \, \frac{(z+1)^3}{4} \, \frac{z^2+2v^{(k)}z+1}{2(v^{(k)}+1)} \,
\\\\
& \hspace{1.4cm} \cdot \, \frac{-\left(v^{(k)}+2 \left(\sqrt{\frac{v^{(k)}+1}{2}}+1 \right) \right)z^2 +
2 \left((v^{(k)} + 1)^2 + 2(v^{(k)}+1) \sqrt{\frac{v^{(k)}+1}{2}} + 1 \right)z
- \left(v^{(k)}+ 2 \left(  \sqrt{\frac{v^{(k)}+1}{2}}+1 \right) \right)}
{4v^{(k)} \sqrt{\frac{v^{(k)}+1}{2}} \left(\sqrt{\frac{v^{(k)}+1}{2}} + 1 \right)}
\end{array}
\end{equation}

\bigskip
\noindent
depending on the parameter $v^{(k)}$ updated at each step through recurrence (\ref{rec_par}).\\
By applying Theorem \ref{theo:ERreproduction} referring to the function space $V_{T,\bgamma}=\{1,\ x, \ e^{tx},\ e^{-tx}\}$ corresponding to the values $\theta_1=0,\ \theta_2=t,\theta_3=-t$ and  $\tau_1=2,\ \tau_2=1,\ \tau_3=1$, with some algebra we easily see that for $$z^{(k)}_1=1,\quad z^{(k)}_2=e^{t/2^{k+1}},\quad z^{(k)}_3=e^{-t/2^{k+1}},\quad t\in \{0,s,\mathrm{i}s\,|\, s>0\}$$ the  conditions
$$
a_3^{(k)}(-z^{(k)}_\ell) =0, \quad a_3^{(k)}(z^{(k)}_\ell)=2(z^{(k)}_\ell)^{-\frac{1}{2}}, \quad \ell=1,\, 2,\, 3,
$$
as well as the conditions
$$
\frac{d}{dz}a_3^{(k)}(-z_1^{(k)})=0, \quad \frac{d}{dz}a_3^{(k)}(z_1^{(k)})=-1
$$
are satisfied so that the scheme reproduces the space $V_{T,\bgamma}$ with respect to the dual parametrization $p=-\frac{1}{2}$.

Since the non-stationary subdivision scheme (\ref{newscheme3}) is asymptotically equivalent \cite{DL95} to the stationary $C^2$ scheme with symbol
$$\frac{(z + 1)^5}{2^4} \left(-\frac{5}{8}z^2+\frac{9}{4}z-\frac{5}{8}\right)=\frac{(z + 1)^2}{2^2} (z+1)\left(-\frac{5}{32}z^4+\frac{1}{4}z^3+\frac{13}{16}z^2+\frac14 z-\frac{5}{32}\right),$$ it is $C^2$ as well.}
\end{example}

\medskip
A last example is given, again, to show the potentialities of Theorem \ref{theo:ERreproduction} in the design of novel non-stationary subdivision schemes with prescribed reproduction properties. In fact
 we now introduce a non-stationary approximating scheme able to reproduce a more complicated function space like $V_{T,\bgamma}=\{1, \ x, \ e^{tx}, \ e^{-tx}, \ x e^{tx}, \ x e^{-tx} \}$. Note that a primal (interpolatory) scheme reproducing this space of functions has been recently presented in \cite{R09}. On the contrary, we here focus our attention on a dual (approximating) scheme reproducing the same space.
\medskip

\begin{example} {Dual scheme reproducing $V_{T,\bgamma}=\{1, \ x, \ e^{tx}, \ e^{-tx}, \ x e^{tx}, \ x e^{-tx} \}$.}

{\rm
\bigskip
\noindent
Let us consider the $k$-level symbol

\small
$$
a_4^{(k)}(z)=z^{-6} \, \frac{(z+1)^3}{4} \, \frac{(z^2+2v^{(k)}z+1)^2}{2(v^{(k)}+1)} \,
\frac{\beta z^4 + \delta z^3 + \epsilon z^2 + \delta z + \beta}{64 (v^{(k)})^3 (v^{(k)}+1) \left(v^{(k)}+1+\frac{1}{2}(v^{(k)}+3) \sqrt{\frac{v^{(k)}+1}{2}} \right)}
$$
\begin{equation}\label{newscheme4}
\end{equation}

\normalsize
with

\small
\begin{displaymath}
\begin{array}{l}
\beta=4 \left((v^{(k)})^2+5v^{(k)}+2 \right) \sqrt{\frac{v^{(k)}+1}{2}}+8(v^{(k)})^2+17v^{(k)}+6,
\\\\
\delta=-4 \sqrt{\frac{v^{(k)}+1}{2}} \left[2 \left(2(v^{(k)})^3+10(v^{(k)})^2+9v^{(k)}+2 \right)+\sqrt{\frac{v^{(k)}+1}{2}} \left(16(v^{(k)})^2+23v^{(k)}+6 \right) \right],
\\\\
\epsilon=8(v^{(k)}+1) \sqrt{\frac{v^{(k)}+1}{2}} \left(2(v^{(k)})^3 + 8(v^{(k)})^2 + 11v^{(k)} + 2 \right)
\\\\
\hspace{0.3cm}+2 \left(16 (v^{(k)})^4 + 48 (v^{(k)})^3 + 70(v^{(k)})^2 + 41v^{(k)} + 6 \right)
\end{array}
\end{displaymath}

\normalsize
\noindent
and $v^{(k)}$ defined as in (\ref{rec_par}).\\
To verify that the non-stationary subdivision scheme with symbols $\{a_4^{(k)}(z), \ k \geq 0\}$ reproduces functions from the space $V_{T,\bgamma}=\{1, \ x, \ e^{tx}, \ e^{-tx}, \ x e^{tx}, \ x e^{-tx} \}$ corresponding to the values $\theta_1=0,\ \theta_2=t,\theta_3=-t$ and  $\tau_1=2,\ \tau_2=2,\ \tau_3=2$, we need to check that, for
$$z^{(k)}_1=1,\quad z^{(k)}_2=e^{t/2^{k+1}},\quad z^{(k)}_3=e^{-t/2^{k+1}},\quad t\in \{0,s,\mathrm{i}s\,|\, s>0\},$$ it results
$$
a_4^{(k)}(-z^{(k)}_\ell) =0, \quad a_4^{(k)}(z^{(k)}_\ell)=2(z^{(k)}_\ell)^{p}, \quad \ell=1,\, 2,\, 3,
$$
and
$$
\frac{d}{dz}a_4^{(k)}(-z_\ell^{(k)})=0, \quad \frac{d}{dz}a_4^{(k)}(z_\ell^{(k)})=2p \, (z_\ell^{(k)})^{p-1}, \quad \ell=1,\, 2,\, 3.
$$
By a few computations it is easy to see that
this is exactly the case with respect to the parametrization $p=-\frac{1}{2}$.

We finally observe that the non-stationary subdivision scheme (\ref{newscheme4}) is of class $C^{3}$ since asymptotically equivalent \cite{DL95} to the stationary $C^3$ scheme with symbol

\begin{displaymath}
\begin{array}{ll}
& \frac{(z + 1)^7}{2^6} \left(\frac{63}{128}z^4-\frac{91}{32}z^3 + \frac{365}{64}z^2 -\frac{91}{32}z+\frac{63}{128} \right)=\\ \\
& \frac{(z + 1)^3}{2^3} (z+1)\left(\frac{63}{1024}z^7 -\frac{175}{1024}z^6 - \frac{173}{1024}z^5 + \frac{797}{1024}z^4 +\frac{797}{1024}z^3 -\frac{173}{1024}z^2 -\frac{175}{1024}z+\frac{63}{1024} \right).
\end{array}
\end{displaymath}
}
\end{example}

\begin{figure}[h!]
\centering
{\includegraphics[trim= 8mm -2mm 8mm 4mm, clip, width=4.5cm]{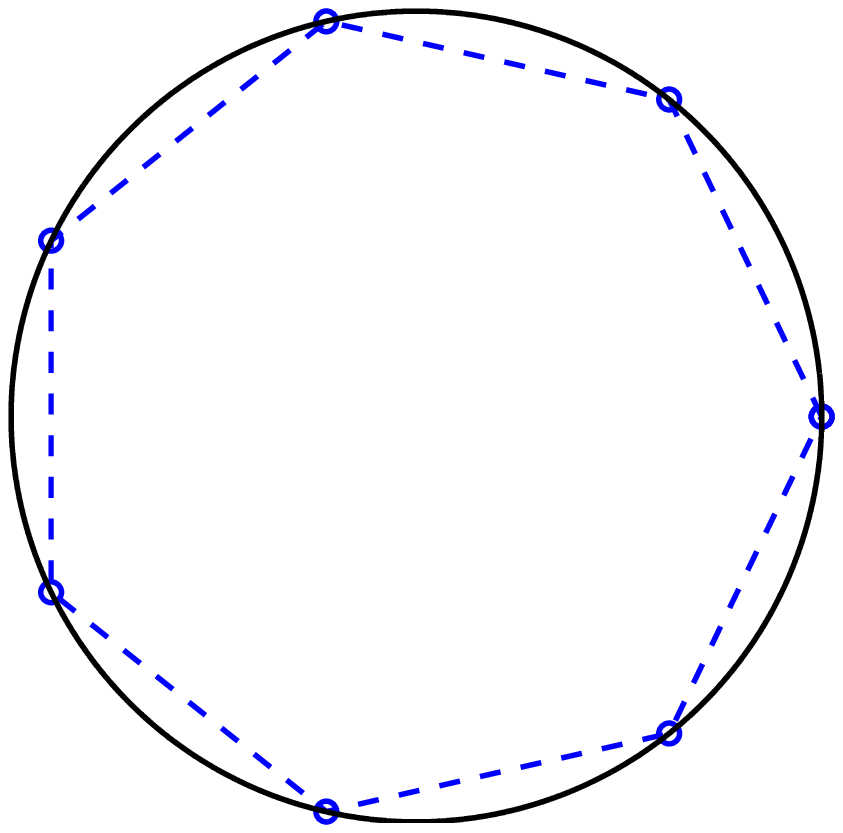}}\hspace{0.2cm}
{\includegraphics[trim= 8mm -2mm 8mm 4mm, clip, width=4.5cm]{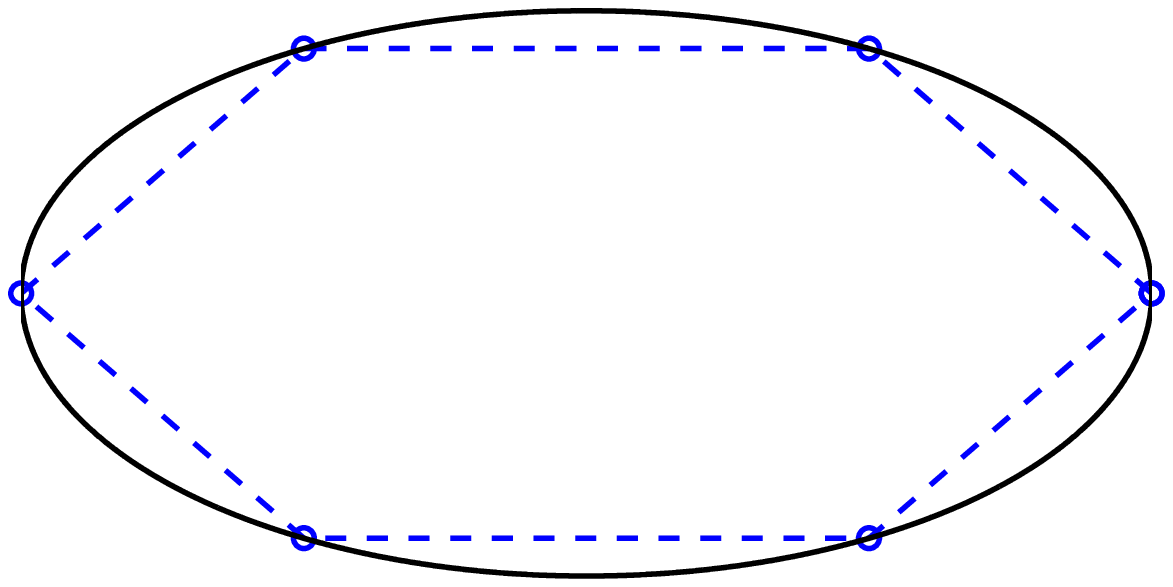}}\hspace{0.2cm}\\
{\includegraphics[trim= 15mm 1mm 8mm 1mm, clip, width=4.2cm]{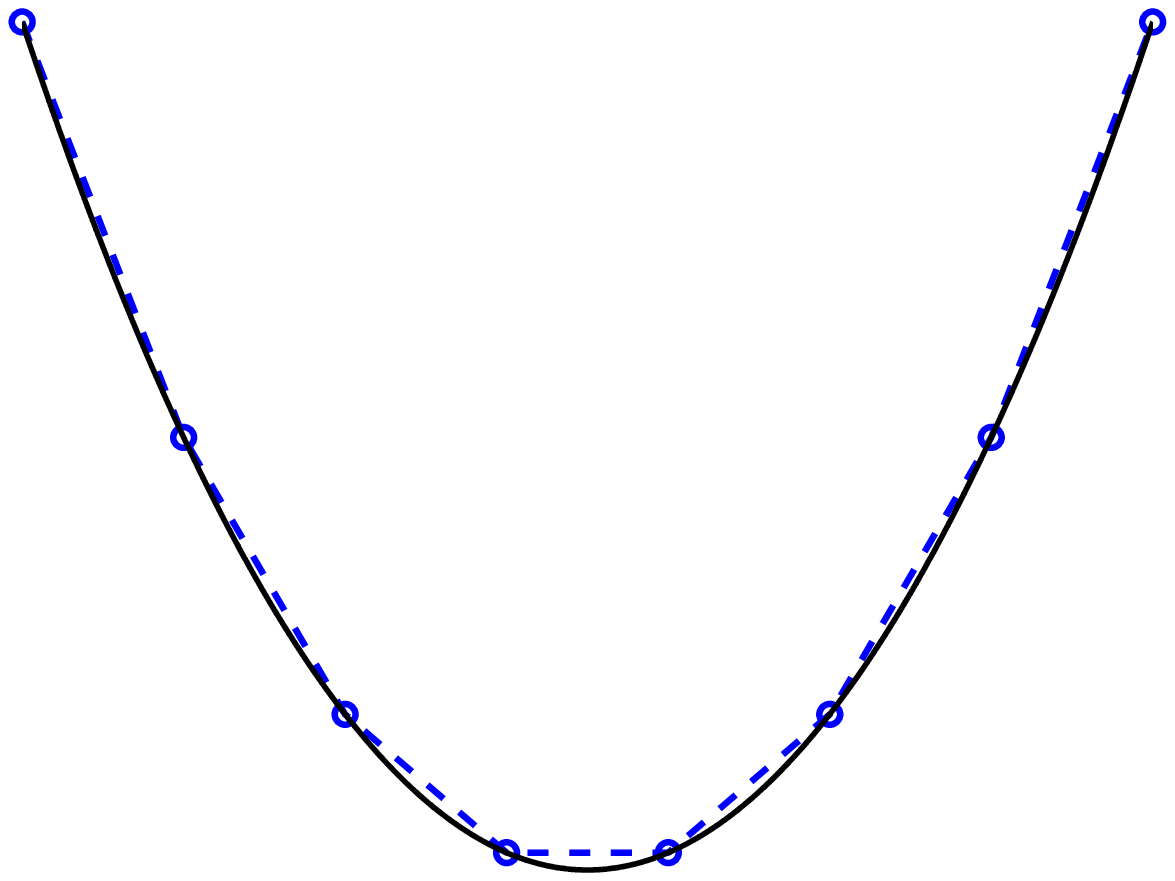}}\hspace{0.2cm}
{\includegraphics[trim= 20mm 2mm 20mm 2mm, clip, width=3.8cm]{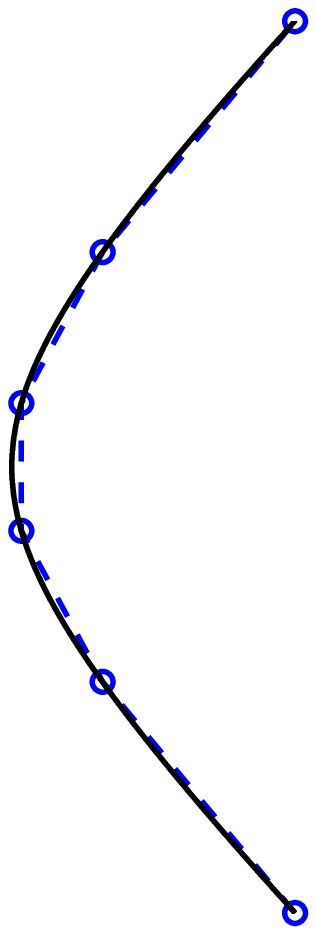}}
\caption{Reproduction of conic sections by the subdivision algorithms with $k$-level symbols in (\ref{newscheme2}), (\ref{newscheme3}) and (\ref{newscheme4}).
The dashed line is the polyline connecting the starting points uniformly sampled from the conic section with constant spacing. The starting parameter $v^{(-1)}$ is chosen equal to $\cos(2\pi/7)$, $\frac{1}{2}$, $1$, $\cosh(3/5)$ for the circle, ellipse, parabola and hyperbola, respectively.}
\label{fig1a}
\end{figure}

\begin{figure}[h!]
\centering
\hspace{-0.5cm}
{\includegraphics[trim= 5mm 5mm 5mm 5mm, clip, width=4.7cm]{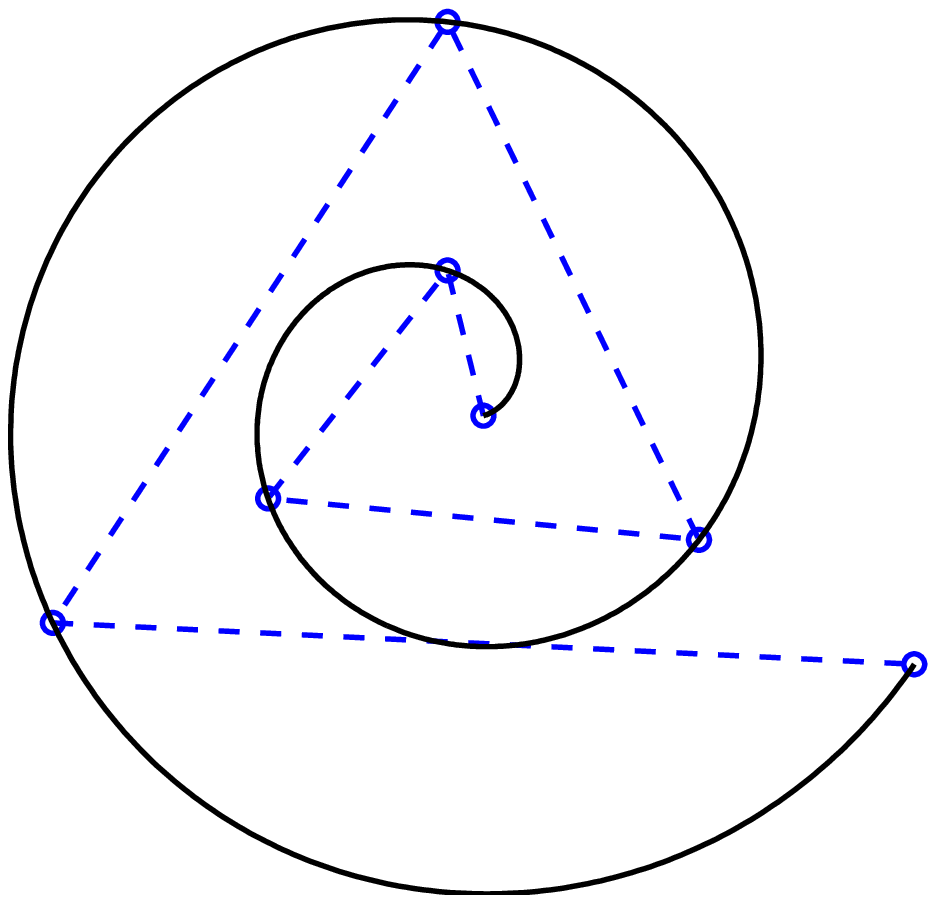}}\hspace{-0.1cm}
{\includegraphics[trim= 5mm 5mm 5mm 5mm, clip, width=4.7cm]{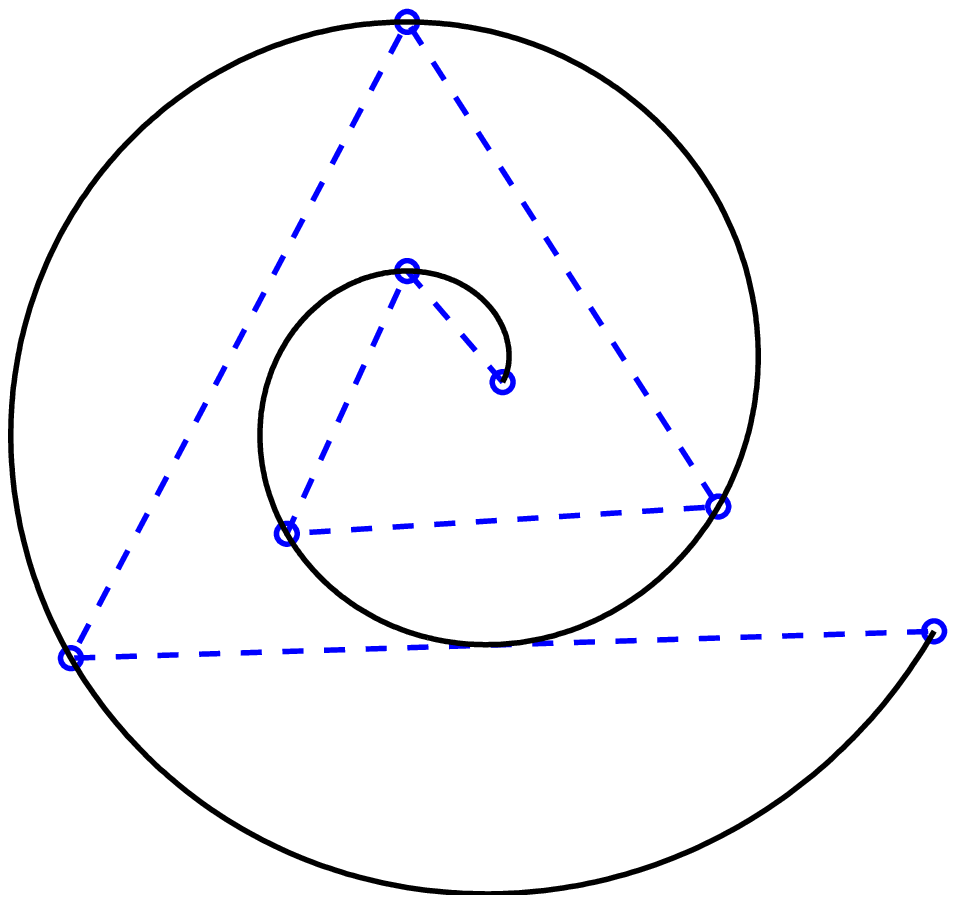}}
\caption{
Reproduction of 2D spirals by the subdivision algorithm with $k$-level symbol in (\ref{newscheme4}).
The dashed line is the polyline connecting the starting points uniformly sampled from the archimedean spiral (left) and the  circle involute (right) with constant spacing. In both cases the starting parameter is $v^{(-1)}=-\frac{1}{2}$. }
\label{fig2a}
\end{figure}

\begin{figure}[h!]
\centering
{\includegraphics[trim= 18mm 4mm 18mm 4mm, clip, width=4.8cm]{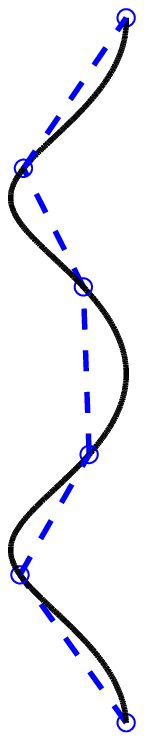}}\hspace{-0.5cm}
{\includegraphics[trim= 6mm 4mm 6mm 4mm, clip, width=4.8cm]{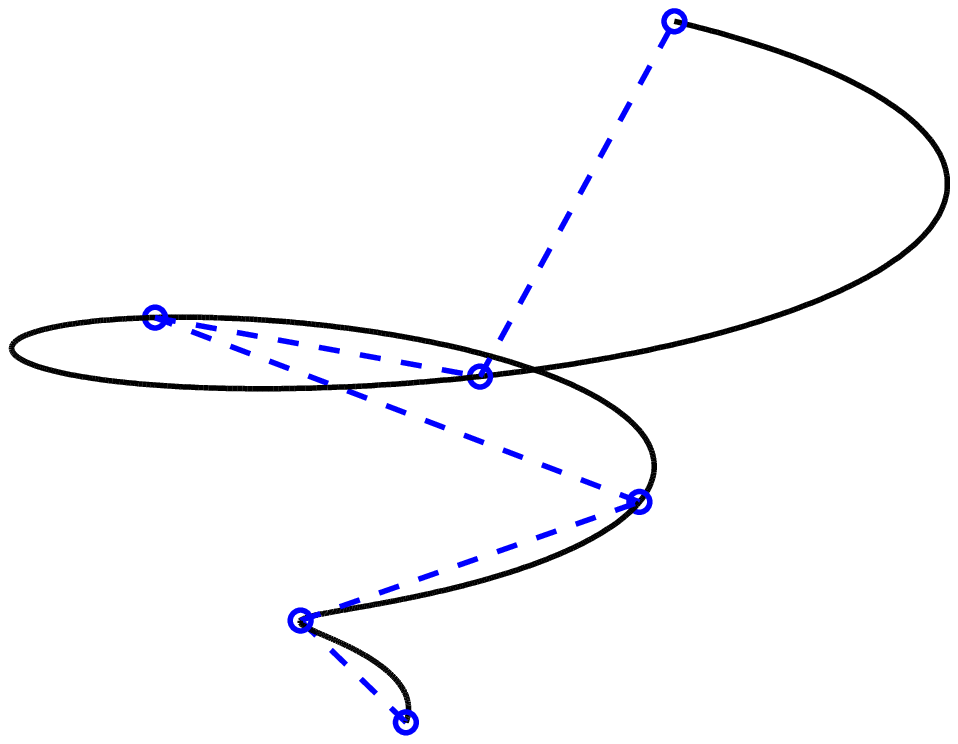}}
\caption{Reproduction of 3D spirals by the subdivision algorithms with $k$-level symbols in (\ref{newscheme2})-(\ref{newscheme4}) (left) and (\ref{newscheme4}) (right).
The dashed line is the polyline connecting the starting points uniformly sampled from the helix (left) and the conical spiral (right) with constant spacing. In both cases the starting parameter is $v^{(-1)}=\cos\left(\frac{4}{5}\pi\right)$.}
\label{fig3a}
\end{figure}

\bigskip
\noindent
Application examples of the proposed schemes to data uniformly sampled from curves of interest in graphical and engineering applications are given in Figures \ref{fig1a}-\ref{fig3a}.
Figure \ref{fig1a} shows that all the three E-reproducing schemes can exactly represent conic sections if starting from points uniformly sampled with constant spacing $\sigma>0$ and choosing $v^{(-1)}=\frac{1}{2}(e^{t \sigma}+e^{-t \sigma})$.
Differently, Figure \ref{fig2a} focuses on a reproduction capability that characterizes the last presented non-stationary scheme.
It consists in the exact reproduction of 2D spirals, such as the archimedean spiral and the circle involute, which are the curves employed in most gear-tooth profiles.
Last, in Figure \ref{fig3a} we extend the reproduction of spirals to the 3D case, by considering the helix and the conical spiral. While the helix can be reconstructed by all the three E-reproducing schemes, the conical spiral can be obtained only by the scheme with symbol in (\ref{newscheme4}).

\end{document}